\documentclass[12pt,a4paper,reqno]{amsart}
\usepackage{amsmath}
\usepackage{amsfonts}
\usepackage{amssymb}
\numberwithin{equation}{section}

\theoremstyle{plain}

\newtheorem{theorem}[subsection]{Theorem}
\newtheorem{proposition}[subsection]{Proposition}
\newtheorem{lemma}[subsection]{Lemma}
\newtheorem{corollary}[subsection]{Corollary}
\newtheorem{conjecture}[subsection]{Conjecture}

\theoremstyle{definition}

\newtheorem{definition}[subsection]{Definition}

\theoremstyle{remark}

\newtheorem*{remark}{Remark}

\renewcommand{\leq}{\leqslant}
\renewcommand{\geq}{\geqslant}

\newsavebox{\proofbox}
\savebox{\proofbox}{\begin{picture}(7,7)%
  \put(0,0){\framebox(7,7){}}\end{picture}}

\def\proof{\noindent\textit{Proof. }}
\def\endproof{\hfill{\usebox{\proofbox}}}

\def\N{\mathbb{N}}
\def\F{\mathbb{F}}
\def\Z{\mathbb{Z}}
\def\R{\mathbb{R}}

\def\eps{\varepsilon}
\def\tildeo{\widetilde{o}}

\newcommand\Spec{\operatorname{Spec}}
\newcommand\Dbl{\operatorname{Dbl}}

\newcommand\Sym{\operatorname{Sym}}

\begin{document}

\title[Freiman and Balog-Szemer\'edi-Gowers theorems]{A note on the Freiman and Balog-Szemer\'edi-Gowers theorems in finite fields}

\author{Ben Green}
\address{Centre for Mathematical Sciences, Wilberforce Road, Cambridge CB3 0WA, England.}
\email{b.j.green@dpmms.cam.ac.uk}

\author{Terence Tao}
\address{Department of Mathematics, UCLA, Los Angeles CA 90095-1555, USA.
}
\email{tao@math.ucla.edu}

\thanks{The first author is a Clay Research Fellow, and is pleased to acknowledge the support of the Clay Mathematics Institute. The second author is supported by a grant from the MacArthur Foundation.}

\begin{abstract} 
We obtain quantitative versions of the Balog-Szemer\'edi-Gowers and Freiman theorems in the model case of a finite field geometry $\F_2^n$, improving the previously known bounds in such theorems.  For instance, if $A \subseteq \F_2^n$ is such that $|A+A| \leq K|A|$ (thus $A$ has small additive doubling), we show that there exists an affine subspace $H$ of $\F_2^n$ of cardinality $|H| \gg K^{-O(\sqrt{K})} |A|$ such that $|A \cap H| \geq \frac{1}{2K} |H|$. Under the assumption that $A$ contains at least $|A|^3/K$ quadruples with $a_1 + a_2 + a_3 + a_4 = 0$ we obtain a similar result, albeit with the slightly weaker condition $|H| \gg K^{-O(K)}|A|$.
\end{abstract}

\maketitle

\section{Introduction}

In this paper we will be working with the group $\F_2^n$, the vector space of dimension $n$ over the two-element field $\F_2$. This serves as a convenient model setting in which to do additive combinatorics. 

If $A,B \subseteq \F_2^n$ are two sets then we define their sumset $A + B$ to be the set of all pairwise sums $a+b$ with $a \in A, b \in B$. A fundamental result concerning sumsets is the following theorem of Ruzsa \cite{ruzsa-group}:

\begin{theorem}[Ruzsa's analogue of Freiman's theorem]\label{ruzsa-theorem}
Let $K \geq 1$ be an integer, and suppose that $A \subseteq \F_2^n$ is a set with $|A + A| \leq K|A|$. Then $A$ is contained in a subspace $H \leq \F_2^n$ with $|H| \leq F(K)|A|$, for some $F(K)$ depending only on $K$.
\end{theorem}

Ruzsa proved this result with $F(K) = K^2 2^{K^4}$. This was then improved by Sanders \cite{sanders-freiman}, who obtained\footnote{As usual we use $X = O(Y)$ or $X \ll Y$ to denote an estimate of the form $X \leq CY$ for some absolute constant $C$.} $F(K) = 2^{O(K^{3/2}\log K)}$. In a recent preprint \cite{green-tao-downsets} the present authors obtain the bound $F(K) = 2^{2K + o(K)}$. This is tight apart from the $o(K)$ term, as can be seen by considering the rather trivial example in which $A$ consists of independent vectors $\{e_1,\dots,e_n\}$ with $n \sim 2K$.

It was already pointed out in \cite{ruzsa-group} that Theorem \ref{ruzsa-theorem} does not quite have the form one would like. Attempting to put \emph{all} of $A$ inside a subspace $H$ is inefficient, as the preceding trivial example illustrates. If one is prepared to place just a portion of $A$ inside a subspace, then conjecturally it is possible to do much better. The following conjecture is attributed in \cite{ruzsa-group} to Marton:

\begin{conjecture}[Polynomial Freiman-Ruzsa Conjecture]\label{pfr}
Suppose that $K \geq 1$ and that $A \subseteq \F_2^n$ has $|A+A| \leq K|A|$. Then there is a subspace $H$ of $\F_2^n$ such that $|H| \ll K^{O(1)} |A|$
and $|A \cap H| \gg K^{-O(1)} |A|$.
\end{conjecture}

It follows from standard covering results in additive combinatorics (see for example \cite[Chapter 2]{tao-vu}) that $A$ is in fact covered by $K^{O(1)}$ translates of $H$. We will not discuss Conjecture \ref{pfr} any further here: the survey \cite{green-pfr} has more information.

It was implicitly observed by Gowers \cite[Chapter 7]{gowers-long-aps} and then by the authors \cite[Chapter 6]{green-tao-u3inverse} that a weak version of Conjecture \ref{pfr} essentially follows from the ideas of Ruzsa \cite{ruzsa-group} on Freiman's theorem in $\Z$. This result was written down (with a somewhat more complicated proof than necessary) in \cite{green-sanders}.

\begin{theorem}[Intersection version of Freiman's theorem in $\F_2^n$]\label{int-frei}
Suppose that $K \geq 1$ and that $A \subseteq \F_2^n$ has $|A+A| \leq K|A|$. Then there is a subspace $H$ of $\F_2^n$ with $|H| \ll K^{O(1)} |A|$ such that $|A \cap H| \gg \exp(-K^{O(1)}) |A|$.
\end{theorem}

This result, of course, is the same as Conjecture \ref{pfr} except that $|A \cap H|/|A|$ may now be exponentially small in $K$. Our first main aim in this paper is to obtain a rather precise version of this result. 

\begin{theorem}[First Main Theorem]\label{main3} Let $K \geq 1$, and let $A,B \subseteq \F_2^n$ be such that $|A+B| \leq K |A|^{1/2} |B|^{1/2}$. Then there exists a linear subspace $H \subseteq \F_2^n$ with $|H| \gg K^{-O(\sqrt{K})}|A|$ and $x,y \in \F_2^n$ such that
\[ |A \cap (x + H)|^{1/2} |B \cap (y + H)|^{1/2} \geq \frac{1}{2K}|H|.\]
\end{theorem}

\begin{remark} We note that the small doubling condition implies that $K^{-2}|A| \leq |B| \leq K^2|A|$, and so we also have $|H| \gg K^{-O(\sqrt{K})}|B|$.\end{remark}

The idea of working with two sets $A,B$ instead of one turns out (for inductive reasons) to be essential to the argument; this idea was suggested to us by Tom Sanders.  It has the following easy corollary:

\begin{corollary}\label{c}  Let $A \subset \F_2^n$ be such that $|A+A| \leq K |A|$ for some $K \geq 1$.  Then there exists a subspace $H \subseteq \F_2^n$ with $|H| \gg K^{-O(\sqrt{K})}|A|$ and $x \in \F_2^n$ such that $|A \cap (x + H)| \geq \frac{1}{2K} |H|$.\endproof
\end{corollary}

Results of Freiman type are particularly powerful when applied in conjunction with statements of Balog-Szemer\'edi type. These results state that if a set $A$ has some weak ``statistical'' additive structure (usually $A$ is assumed to have many \emph{additive quadruples}, i.e. quadruples $(a_1,a_2,a_3,a_4) \in A^4$ with $a_1 + a_2 = a_3 + a_4$) then there is a large subset $A' \subseteq A$ for which $|A' + A'|$ is small.
Such an application was first made in \cite{balog}. A major advance was made by Gowers \cite{gowers-4}, who proved a result of Balog-Szemer\'edi type with good quantitative bounds.

\begin{theorem}[Balog-Szemer\'edi-Gowers theorem]\label{bsg}
Let $G$ be any abelian group, and suppose that $A \subseteq G$ is a finite set with at least $|A|^3/K$ additive quadruples. Then there is $A' \subseteq A$ with $|A'| \gg K^{-O(1)}|A|$ such that $|A' + A'| \ll K^{O(1)}|A'|$. 
\end{theorem}

This theorem is now proven in many places in the literature, see e.g. \cite[Theorem 2.29]{tao-vu}.
In combination with Theorem \ref{int-frei}, the Balog-Szemer\'edi-Gowers theorem implies the following result.

\begin{theorem}[Balog-Szemer\'edi-Gowers-Freiman theorem in $\F_2^n$]
Suppose that $A \subseteq \F_2^n$ is a set with at least $|A|^3/K$ additive quadruples. Then there is a subspace $H$ of $\F_2^n$ with $|H| \gg \exp(-K^{O(1)})|A|$ and an $x \in \F_2^n$ such that $|A \cap (x + H)| \gg \exp(-K^{O(1)})|H|$.
\end{theorem}

Our second main aim in this paper is to prove a more precise version of this latter result. We shall need some notation:

\begin{definition}[Normalised energy]
Given any non-empty sets $A_1,A_2,A_3,A_4 \subseteq \F_2^n$, define the \emph{normalised energy}
\begin{equation}\label{energy-def}
\omega(A_1,A_2,A_3,A_4) := \frac{|\{ (a_1,a_2,a_3,a_4) \in A_1 \times A_2 \times A_3 \times A_4 : a_1 + a_2 + a_3 + a_4 = 0 \}|}{(|A_1| |A_2| |A_3| |A_4|)^{3/4}}.
\end{equation}
\end{definition}

Thus for instance the statement that $A$ has at least $|A|^3/K$ additive quadruples is equivalent to the assertion that
$$ \omega(A,A,A,A) \geq 1/K$$
(note that in the characteristic $2$ group $\F_2^n$, there is no distinction between addition and subtraction). Noting that the number of quadruples with $a_1 + a_2 + a_3 + a_4 = 0$ is bounded by $\prod_{i \neq j} |A_i|$ for any $j$, we see by taking products over $j = 1,\dots,4$ that $0 \leq \omega(A_1,A_2,A_3,A_4) \leq 1$.  A simple application of Cauchy-Schwarz also gives the inequality
\begin{equation}\label{omegab}
\omega(A,B,A,B) \geq \frac{|A|^{1/2} |B|^{1/2}}{|A+B|}.
\end{equation}

\begin{theorem}[Second Main Theorem]\label{main4}
Let $A_1,A_2,A_3,A_4 \subseteq \F_2^n$ be non-empty sets such that
\[ \omega(A_1,A_2,A_3,A_4) \geq \frac{1}{K}\]
for some $K \geq 1$.  Then there exists a linear subspace $H \subseteq \F_2^n$ with $|H| \gg K^{-O(K)}|A_i|$ for $i = 1,2,3,4$ and $x_1,x_2,x_3,x_4 \in \F_2^n$ such that
 \[ \prod_{i=1}^4 |A_i \cap (x_i + H)|^{1/4} \geq \frac{1}{2K}|H|.\]
\end{theorem}

The proof of this result is closely related to that of Theorem \ref{main3}. More importantly it is \emph{direct} in the sense that the Balog-Szemer\'edi-Gowers theorem is not required. It therefore demonstrates that the Balog-Szemer\'edi-Gowers and Freiman theorems, which are traditionally proven using very different methods, can instead be treated in a unified manner. Note however that we do not recover the polynomial-type bounds in Theorem \ref{bsg} by our methods. Once again there is a simple corollary when just one set is involved:

\begin{corollary}  Let $A \subseteq \F_2^n$ have at least $|A|^3/K$ additive quadruples for some 
$K \geq 1$.  Then there exists an affine subspace $H \subseteq \F_2^n$ with
$|H| \gg K^{-O(K)} |A|$ and $|A \cap H| \geq \frac{1}{2K} |H|$.\endproof
\end{corollary}

Our arguments used to prove Theorem \ref{main3} and Theorem \ref{main4} are completely self-contained, and can be briefly summarised as follows.
First, in \S \ref{well-sec}, we use Fourier-analytic methods to establish a special case of Theorem \ref{main4} (or Theorem \ref{main3}) when the four sets $A_1,A_2,A_3,A_4$ are ``coherently flat'', in the sense that the Fourier coefficients are either simultaneously large or simultaneously small. See in particular Proposition \ref{well}.  Then, in \S \ref{freiman-sec}, we run an energy increment argument to reduce to this coherently flat case in Theorem \ref{main3}. In \S \ref{balog-sec} we perform a similar argument to prove Theorem \ref{main4}.  

\section{The coherently flat case}\label{well-sec}

Our first task is to prove a variant of Theorems \ref{main3} and \ref{main4} in the case where we have four sets $A_1,A_2,A_3,A_4$ which are ``coherently flat'', which basically corresponds to the regime in which Fourier analysis tools are effective.  

\begin{definition}[Fourier transform] If $f: \F_2^n \to \R$, we define the Fourier transform $\hat f: \F_2^n \to \R$ by the formula
\[ \hat f(\xi) := \frac{1}{2^n} \sum_{x \in \F_2^n} f(x)(-1)^{\xi \cdot x} \]
where $(x_1,\ldots,x_n) \cdot (\xi_1,\ldots,\xi_n) := x_1 \xi_1 + \ldots + x_n \xi_n$.
\end{definition}

Now we use a definition from \cite{tao-vu}.

\begin{definition}[Spectrum]\label{spec-def}
If $A \subseteq \F_2^n$ is non-empty and $0 < \alpha \leq 1$, define the \emph{$\alpha$-spectrum} $\Spec_\alpha(A)$ to be the set of all frequencies $\xi \in \F_2^n$ such that
\[ |\hat 1_A(\xi)| \geq \alpha |A|/2^n.\]
\end{definition}

The spectrum $\Spec_\alpha(A)$ can be viewed as collecting the directions in which $A$ is significantly biased; indeed to say that $\xi \in \Spec_\alpha(A)$ is equivalent to the assertion that the proportion of $A$ in the subspace $\{ x \in \F_2^n: \xi \cdot x = 0 \}$ is either greater than $\frac{1+\alpha}{2}$, or less than $\frac{1-\alpha}{2}$.

Observe that if $A$ is an affine subspace, then $\hat 1_A(\xi)$ is either equal to $|A|/2^n$ or zero; in other words, given any $\xi$, either $\xi \in \Spec_1(A)$ or $\xi \not \in \Spec_\delta(A)$ for all $\delta > 0$.  Slightly more generally, we see that if $A_1,A_2,A_3,A_4$ are cosets of the same subspace, then for any $\xi$, either $\xi \in \Spec_1(A_1) \cap \ldots \cap \Spec_1(A_4)$, or
$\xi\not \in \Spec_\delta(A_1) \cup \ldots \cup \Spec_\delta(A_4)$ for all $\delta > 0$.  This motivates the following definition of a quadruple of sets which ``resemble'' four cosets of the same subspace in some Fourier sense.

\begin{definition}[Coherently flat quadruples]
Suppose that $A_1,A_2,A_3,A_4 \subseteq \F_2^n$ are non-empty and that $\delta \in (0,\frac{1}{2})$ is a small parameter. We say that the quadruple $(A_1,A_2,A_3,A_4)$ is \emph{coherently $\delta$-flat} if, for each $\xi \in \F_2^n$, one of the following is true:
\begin{enumerate}
\item \textup{($\xi$ orthogonal to all $A_i$)} $\xi \in \Spec_{9/10}(A_i)$ for all $i=1,2,3,4$;
\item \textup{($\xi$ non-orthogonal to all $A_i$)} $\xi \notin \Spec_{\delta}(A_i)$ for all $i=1,2,3,4$.
\end{enumerate}
\end{definition}

We observe that we may translate one or more of the $A_i$ by an arbitrary $x_i \in \F_2^n$ without affecting the coherent flatness property.

The main result of this section is as follows.

\begin{proposition}[Freiman-type theorem, coherently flat case]\label{well}  Let $K \geq 1$.  
Suppose $(A_1,A_2,A_3,A_4)$ is a coherently $\frac{1}{\sqrt{2K}}$-flat quadruple whose energy satisfies the lower bound $\omega(A_1,A_2,A_3,A_4) \geq 1/K$.
Then there is a linear subspace $H \subseteq \F_2^n$ with 
\begin{equation}\label{hlower}
|H| \geq \frac{4}{5} \prod_{i=1}^4 |A_i|^{1/4},
\end{equation}
together with $x_1,x_2,x_3,x_4 \in \F_2^n$
 such that 
\begin{equation}\label{alower}
\prod_{i=1}^4 |A_i \cap (x_i + H)|^{1/4} \geq \frac{1}{2K}|H|.
\end{equation}
\end{proposition}

\begin{proof}  Set
$$ \Lambda := \Spec_{9/10}(A_1) \cap \dots \cap \Spec_{9/10}(A_4).$$
We claim that $\Lambda$ is a linear subspace of $\F_2^n$.  Indeed, if this were not the case then we could find $\xi, \xi' \in \Lambda$ such that $\xi + \xi' \not \in \Lambda$.  Without loss of generality we may assume that $\xi + \xi' \not \in \Spec_{9/10}(A_1)$.  Then by the coherently $\frac{1}{\sqrt{2K}}$-flat hypothesis
we see that $\xi + \xi' \not \in \Spec_{1/\sqrt{2K}}(A_1)$.  On the other hand, we have $\xi, \xi' \in \Spec_{9/10}(A_1)$, which by the triangle inequality (cf. \cite[Lemma 4.37]{tao-vu}) implies $\xi+\xi' \in \Spec_{8/10}(A_1)$. This is a contradiction, and so $\Lambda$ is indeed a linear subspace.  

Using the Fourier transform (see e.g. \cite[Lemma 4.9]{tao-vu}) we observe the identity
\begin{align*}
 \sum_{\xi \in \F_2^n}
\prod_{i=1}^4 \widehat{1}_{A_i}(\xi) 
&= \frac{\prod_{i=1}^4 |A_i|^{3/4}}{2^{3n}} \omega(A_1,A_2,A_3,A_4) \\
&\geq \frac{\prod_{i=1}^4 |A_i|^{3/4}}{2^{3n} K}.
\end{align*}
On the other hand the coherently $\frac{1}{\sqrt{2K}}$-flat hypothesis, H\"older's inequality, and Plancherel imply that
\begin{align*}
\big|\sum_{\xi \in \F_2^n \backslash \Lambda}
\prod_{i=1}^4 \widehat{1}_{A_i}(\xi)\big|
&\leq \sum_{\xi \in \F_2^n \backslash \Lambda}
\prod_{i=1}^4 \left( \frac{|A_i|}{2^n \sqrt{2K}} \right)^{1/2}
|\widehat{1}_{A_i}(\xi)|^{1/2} \\
&\leq 
\prod_{i=1}^4 \left( \frac{|A_i|}{2^n \sqrt{2K}} \right)^{1/2}
(\sum_{\xi \in \F_2^n} |\widehat{1}_{A_i}(\xi)|^2)^{1/4} \\
&= \prod_{i=1}^4 \left( \frac{|A_i|}{2^n \sqrt{2K}} \right)^{1/2} \left( \frac{|A_i|}{2^n} \right)^{1/4} \\
&= \frac{1}{2}\frac{\prod_{i=1}^4 |A_i|^{3/4}}{2^{3n} K}.
\end{align*}
Subtracting the two estimates, we conclude that
\begin{equation}\label{xih}
 \sum_{\xi \in \Lambda}
\prod_{i=1}^4 \widehat{1}_{A_i}(\xi) \geq \frac{1}{2}\frac{\prod_{i=1}^4 |A_i|^{3/4}}{2^{3n} K}.
\end{equation}
Write $H$ for the orthogonal complement of $\Lambda$, that is to say
\[H := \Lambda^\perp := \{ x \in \F_2^n: x \cdot \xi = 0 \hbox{ for all } \xi \in \Lambda \}.\]
The left-hand side of \eqref{xih} can be rewritten as
\[ \frac{1}{2^{3n} |H|} | \{ (a_1,a_2,a_3,a_4) \in A_1 \times A_2 \times A_3 \times A_4:
a_1 + a_2 + a_3 + a_4 \in H \} |.\]
This is bounded above by
\[ \frac{1}{2^{3n} |H|} |A_2| |A_3| |A_4| \sup_{x_1 \in \F_2^n} |A_1 \cap (x_1 + H)|\]
It follows from \eqref{xih} that there exists $x_1 \in \F_2^n$ such that
\[ |A_1 \cap (x_1 + H)| |A_2| |A_3| |A_4| \geq \frac{|H|}{2K} (|A_1| \ldots |A_4|)^{3/4}.\]
We may similarly find $x_2,x_3,x_4$ such that similar inequalities hold with the $A_i$ permuted. Taking geometric means of these four estimates and rearranging we obtain \eqref{alower}.

It remains to prove \eqref{hlower}.  Observe from Plancherel's theorem and the definition of $\Lambda$ that for any $i=1,2,3,4$ we have
\[ \frac{|A_i|}{2^n} = \sum_{\xi \in \F_2^n} |\hat 1_{A_i}(\xi)|^2
\geq \sum_{\xi \in \Lambda} \left( \frac{9}{10} \frac{|A_i|}{2^n} \right )^2.\]
It follows that
\[ |\Lambda| \leq \frac{5 \times 2^n}{4 |A_i|},\]
and so upon taking geometric means we conclude that
\[ |\Lambda| \leq \frac{5 \times 2^n}{4 \prod_{i=1}^4 |A_i|^{1/4}}.\]
Since $|H| = 2^n / |\Lambda|$, the claim \eqref{hlower} follows.
\end{proof}

\section{Proof of Theorem \ref{main3}}\label{freiman-sec}

We now prove Theorem \ref{main3}.  We begin by proving an elementary lemma.  
Let $F: [0,1] \times [0,1] \to \R^+$ be the explicit function
$$ F(x,y) := \sqrt{x} ( \sqrt{y} + \sqrt{1-y} ).$$
Thus for instance $F(1/2,1/2) = F(1,1) = F(1,0) = 1$.  

\begin{lemma}[Near-minima of $F$]\label{lem0.3}
Suppose that $\alpha,\beta,\epsilon$ are reals with $0 \leq \alpha,\beta \leq 1$ and $0 < \eps \leq \frac{1}{100}$, and that
\begin{equation}\label{squaw}
F(\alpha,\beta), F(1-\alpha,\beta), F(\beta,\alpha), F(\beta,1-\alpha) \leq 1+\eps.
\end{equation}
Then either both $\alpha$ and $\beta$ are within $4\eps$ of $\frac{1}{2}$, or else they are both within $2\eps$ of $0$ or $1$.
\end{lemma}

\begin{proof} The hypotheses are invariant under interchanging $\alpha$ and $\beta$, or by swapping $\alpha$ to $1-\alpha$ or $\beta$ to $1-\beta$.  Thus without loss of generality we may assume that $\frac{1}{2} \leq \alpha \leq \beta$.
We have 
\[ \alpha + \sqrt{\alpha(1-\alpha)} = F(\alpha,\alpha) \leq F(\beta,\alpha) \leq 1 + \eps;\]
using the inequality $\sqrt{\alpha(1-\alpha)} \geq 2\alpha(1-\alpha)$ and rearranging we obtain
\[ (2\alpha - 1)(1 - \alpha) \leq \eps.\] 
It follows immediately that $\alpha$ is within $2\eps$ of either $\frac{1}{2}$ or $1$.
In the latter case we are done, so suppose that $\alpha$ is within $2\eps$ of $\frac{1}{2}$. We have
\[ F(\beta,\alpha)^2 \leq (1 + \eps)^2 \leq 1+3\eps\]
which expands to
\[ \beta \big(1 + 2\sqrt{\frac{1}{4} - (\frac{1}{2} - \alpha)^2} \big) \leq 1 + 3\eps\] 
and whence
\[ \beta \big(1 + \sqrt{1 - 16\eps^2}\big) \leq 1 + 3\eps.\]
This implies that $\beta$ is within $4\eps$ of $\frac{1}{2}$, and the claim follows.
\end{proof}

As a consequence of this lemma we obtain the following key inductive step required for Theorem \ref{main3}.  It is convenient to introduce the notation $\Dbl(A,B) := \frac{|A+B|}{|A|^{1/2} |B|^{1/2}}$.  Since $|A+B| \geq \max(|A|, |B|)$ we observe that
\begin{equation}\label{distance}
\Dbl(A,B)^{-2} |A| \leq |B| \leq \Dbl(A,B)^2 |A|.
\end{equation}

\begin{lemma}[Non-flatness implies doubling decrement]\label{lem0.4} Suppose that $A,B \subseteq \F_2^n$ are non-empty sets with $\Dbl(A,B) \leq K$ for some $K \geq 1$, and suppose that $(A,B,A,B)$ is not coherently $\frac{1}{\sqrt{2K}}$-flat. Then there are $A' \subseteq A$, $B' \subseteq B$ with $|A'|/|A|, |B'|/|B| \gg K^{-10}$ and 
\begin{equation}\label{star} \Dbl(A',B') \leq K - \frac{1}{100} \sqrt{K}.\end{equation}
\end{lemma}

\begin{proof} By hypothesis, we can find $\xi \in \F_2^n$ such that
\begin{equation}\label{xihigh}
\xi \notin \Spec_{9/10}(A) \cap \Spec_{9/10}(B)
\end{equation}
and
\begin{equation}\label{xilow}
\xi \in \Spec_{1/\sqrt{2K}}(A) \cup \Spec_{1/\sqrt{2K}}(B).
\end{equation}
Observe that $\xi$ must be non-zero.

For $j \in \F_2$ we set $A_j := \{x \in A: x \cdot \xi = j \}$ and $B_j := \{x \in B: x \cdot \xi = j \}$,
and write $\alpha := |A_0|/|A|$ and $\beta := |B_0|/|B|$.   Then
\[ |\hat 1_A(\xi)| = |2\alpha - 1| \frac{|A|}{2^n} \qquad \mbox{and} \qquad|\hat 1_B(\xi)| = |2\beta - 1| \frac{|B|}{2^n}\]
and so \eqref{xihigh} is equivalent to the assertion
$$
|2\alpha - 1| < \frac{9}{10} \qquad \mbox{or} \qquad |2\beta - 1| < \frac{9}{10}
$$
while \eqref{xilow} is equivalent to the assertion
$$
 |2\alpha - 1| \geq \frac{1}{\sqrt{2K}} \qquad \mbox{or} \qquad |2\beta - 1| \geq \frac{1}{\sqrt{2K}}.
$$
Applying Lemma \ref{lem0.3} in the contrapositive we conclude that one of the four quantities 
$F(\alpha,\beta)$, $F(1-\alpha,\beta)$, $F(\beta,\alpha)$, $F(1-\beta,\alpha)$  is greater than $1 + \frac{1}{100\sqrt{K}}$.  By swapping $A$ and $B$, or by shifting either $A$ or $B$ by $\xi$ we may assume without loss of generality that
\[ F(\alpha,\beta) > 1 + \frac{1}{100\sqrt{K}}.\]
In particular we see that $\beta \neq 0,1$ and $\alpha \neq 0$, so that $A_0,B_0,B_1$ are non-empty; a variant of this argument also gives $|A_0| \gg K^{-10} |A|$ and $|B_j| \gg K^{-10} |B|$ for $j \in \F_2$. We can rewrite the above inequality as
$$ \sum_{j \in \F_2} \frac{K}{1 + \frac{1}{100\sqrt{K}}}
|A_0|^{1/2} |B_j|^{1/2} > K |A|^{1/2} |B|^{1/2}.$$
On the other hand, we have
$$ K |A|^{1/2} |B|^{1/2} \geq |A+B| \geq \sum_{j \in \F_2} |A_0 + B_j|,$$
and so there exists $j \in \F_2$ such that
$$ \Dbl(A_0,B_j) \leq \frac{K}{1 + \frac{1}{100\sqrt{K}}}
\leq K - \frac{1}{100} \sqrt{K}.$$
The claim follows.
\end{proof}

We can iterate the above lemma at most $O(\sqrt{K})$ times, noting that the doubling constant cannot drop below $1$, to obtain

\begin{corollary}[Large coherently flat subsets]\label{corflat} Suppose that $A,B \subseteq \F_2^n$ are non-empty sets with $\Dbl(A,B) \leq K$ for some $K \geq 1$.  Then there exists $A' \subseteq A$, $B' \subseteq B$ with $|A'|/|A|, |B'|/|B| \gg K^{-O(\sqrt{K})}$, $\Dbl(A',B') \leq K$, and such that $(A',B',A',B')$ is coherently $\frac{1}{\sqrt{2K}}$-flat.
\end{corollary}

\emph{Proof of Theorem \ref{main3}.}  We apply Corollary \ref{corflat} to extract $A', B'$ with the stated properties.  
Now we observe the identities
$$ \omega(A',B',A',B') = \frac{1}{ (|A'||B'|)^{3/2} } 
\sum_{x \in A'+B'} |\{ (a,b) \in A' \times B': a + b \}|^2$$
and
$$ \sum_{x \in A'+B'} |\{ (a,b) \in A' \times B': a + b \}| = |A'| |B'|.$$
Applying the Cauchy-Schwarz inequality we conclude that
$$ \omega(A',B',A',B') \geq \frac{|A'|^{1/2} |B'|^{1/2}}{|A'+B'|} = \frac{1}{\Dbl(A',B')} \geq \frac{1}{K}.$$
Applying Proposition \ref{well}, we find $H,x_1,x_2,x_3,x_4$ such that
$$
|H| \geq \frac{4}{5} |A'|^{1/2} |B'|^{1/2} \geq K^{-O(\sqrt{K})} |A|^{1/2} |B|^{1/2}$$
and
$$ |A \cap (x_1 + H)|^{1/4} |B \cap (x_2 + H)|^{1/4}
 |A \cap (x_3 + H)|^{1/4} |B \cap (x_4 + H)|^{1/4}
 \geq \frac{1}{2K} |H|.$$
Without loss of generality we may assume that $|A \cap (x_1+H)| \geq |A \cap (x_3+H)|$ and
$|B \cap (x_2+H)| \geq |B \cap (x_4+H)|$, thus
$$ |A \cap (x_1+H)|^{1/2} |B \cap (x_2+H)|^{1/2} \geq \frac{1}{2K} |H|.$$
The claim follows.
\endproof

\section{Proof of Theorem \ref{main4}}\label{balog-sec}

We can adapt the arguments of the previous section to prove Theorem \ref{main4}.  We begin with the
analogue of Lemma \ref{lem0.3}.  We now require the explicit function $G: [0,1]^4 \to \R^+$ defined by
$$ G(\alpha_1,\alpha_2,\alpha_3,\alpha_4) := \sum_{\substack{j_1,j_2,j_3,j_4 \in \F_2 \\ j_1+j_2+j_3+j_4 = 0}} \prod_{i=1}^4 \alpha_{i,j_i}^{3/4},$$
where we write $\alpha_{i,0} := \alpha_i$ and $\alpha_{i,1} := 1-\alpha_i$.

Applying Young's inequality in the form
$$ |f_1 * f_2 * f_3 * f_4(0)| \leq \prod_{i=1}^4 (\sum_{j \in \F_2} |f_i(j)|^{4/3})^{3/4}$$
with $f_i(j) := \alpha_{i,j}^{3/4}$ we see that $G$ is always bounded by $1$.  Equality occurs in this inequality if and only if each of the $f_i$ are concentrated on a single point in $\F_2$, or else are all constant in $\F_2$.  We shall need a robust version of this observation 

\begin{lemma}[Near-maxima of $G$]\label{olympiad} Suppose that $\alpha_1,\alpha_2,\alpha_3,\alpha_4,\epsilon$ are reals and that $0 \leq \alpha_1,\ldots,\alpha_4 \leq 1$ and $0 < \eps \leq \frac{1}{1000}$. Suppose that
$$
G(\alpha_1,\alpha_2,\alpha_3,\alpha_4) \geq 1-\eps.$$
Then either the four quantities $\alpha_i$ are all within $3\sqrt{\eps}$ of $\frac{1}{2}$, or else the four quantities $\min(\alpha_i, 1-\alpha_i)$ are all at most $10\eps$.
\end{lemma}

\begin{proof} Observe that
$$ G(\alpha_1,\alpha_2,\alpha_3,\alpha_4) = 
\frac{1}{2} \prod_{i=1}^4 (\alpha_i^{3/4} + (1 - \alpha_i)^{3/4})
+ \frac{1}{2} \prod_{i=1}^4 (\alpha_i^{3/4} - (1 - \alpha_i)^{3/4}).$$

Applying the arithmetic-geometric mean inequality
\[ x_1 x_2 x_3 x_4 \leq \frac{x_1^4 + x_2^4 + x_3^4 + x_4^4}{4},\]
we obtain
\[ 8 -  \sum_{i=1}^4 \big( \alpha_i^{3/4} + (1 - \alpha_i)^{3/4}\big)^4  - 
\sum_{i=1}^4 \big( \alpha_i^{3/4} - (1 - \alpha_i)^{3/4}\big)^4 \leq 8\eps.\]
After some rearrangement, the left hand side may be seen to equal
\[ \sum_{i=1}^4 \frac{3\alpha_i(1 - \alpha_i)(1 - 2\alpha_i)^2}{1 + 2\alpha_i^{1/2}(1 - \alpha_i)^{1/2}},\]
which is at least
\[ \sum_{i=1}^4 \frac{3}{2}\alpha_i(1 - \alpha_i)(1 - 2\alpha_i)^2.\]
By positivity we thus have
\[ |\alpha_i||1 - \alpha_i||\frac{1}{2} - \alpha_i|^2 \leq 4\eps/3,\]
for $i=1,2,3,4$.  A short back of an envelope computation confirms that each of $\alpha_i$ is either within $3\sqrt{\eps}$ of $\frac{1}{2}$, or else within $10\eps$ of $0$ or $1$.

This is not quite as strong as the statement of the lemma, because different $\alpha_i$ may end up in different sides of the dichotomy.  Suppose for contradiction that this occurs. Without loss of generality that $\alpha_1$ is within $3\sqrt{\eps}$ of $\frac{1}{2}$. Then we have
\[ \frac{1}{2} \prod_{i=1}^4 |\alpha_i^{3/4} - (1 - \alpha_i)^{3/4}| \leq 0.085\] and so
\[ \frac{1}{2} \prod_{i=1}^4 (\alpha_i^{3/4} +(1 - \alpha_i)^{3/4}) \geq \frac{9}{10}.\]
Since $x^{3/4} + (1 - x)^{3/4} \leq 2^{1/4}$ for $0 \leq x \leq 1$, we obtain
\[ \alpha_i^{3/4} + (1 - \alpha_i)^{3/4} \geq \frac{9}{5 \cdot 2^{3/4}} \geq 1.07\]
for $i=2,3,4$. It follows that none of $\alpha_2,\alpha_3,\alpha_4$ is within $10\epsilon$ of $0$ or $1$, and hence all of these must be within $3\sqrt{\eps}$ of $\frac{1}{2}$ as well.
\end{proof}

Now we obtain the analogue of Lemma \ref{lem0.4}.

\begin{lemma}[Non-flatness implies energy increment]\label{lem1.4} Let $K \geq 1$.
Suppose that $A_1$, $A_2,A_3,A_4 \subseteq \F_2^n$ are non-empty sets with $\omega(A_1,A_2,A_3,A_4) \geq 1/K$, and suppose that $(A_1,A_2,A_3,A_4)$ is not coherently $\frac{1}{\sqrt{2K}}$-flat.
Then for each $i=1,2,3,4$ there are sets $A'_i \subseteq A_i$ with $|A'_i|/|A_i| \gg K^{-10}$ such that 
\begin{equation}\label{star2} \omega(A'_1,A'_2,A'_3,A'_4) \geq \frac{1}{K - 10^{-4}}.\end{equation}
\end{lemma}

\begin{proof} By hypothesis, we can find $\xi$ such that
$$
\xi \notin \Spec_{9/10}(A_1) \cap \ldots \cap \Spec_{9/10}(A_4)
$$
and
$$
\xi \in \Spec_{1/\sqrt{2K}}(A_1) \cup \ldots \cup \Spec_{1/\sqrt{2K}}(A_4).
$$
Again $\xi$ must be non-zero.  We then set $A_{i,j} := \{ x\in A_i: \xi \cdot x = j \}$ for $i=1,2,3,4$ and $j \in \F_2$, and set
$\alpha_i := |A_{i,0}|/|A_i|$.  Then we have
\begin{equation}\label{alphamin}
 \min_{i=1,2,3,4} |2\alpha_i - 1| < 9/10
 \end{equation}
and
\begin{equation}\label{alphamax}
 \max_{i=1,2,3,4} |2\alpha_i - 1| \geq 1/\sqrt{2K}.
\end{equation}
Applying Lemma \ref{olympiad} in the contrapositive we conclude that
$$
G(\alpha_1,\alpha_2,\alpha_3,\alpha_4) < 1 - \frac{1}{1000 K},$$
and consequently
$$
\sum_{\substack{j_1,j_2,j_3,j_4 \in \F_2\\ j_1+j_2+j_3+j_4 = 0}} \big(\prod_{i=1}^4 |A_{i,j_i}|^{3/4}
+ \frac{1}{10^4 K} \prod_{i=1}^4 |A_i|^{3/4}\big)
 < (1 - \frac{1}{10^4 K}) \prod_{i=1}^4 |A_i|^{3/4}.$$
Now observe the identity
$$ \omega(A_1,A_2,A_3,A_4) \prod_{i=1}^4 |A_i|^{3/4} =
\sum_{\substack{j_1,j_2,j_3,j_4 \in \F_2\\ j_1+j_2+j_3+j_4 = 0} }
\omega(A_{1,j_1},A_{2,j_2},A_{3,j_3},A_{4,j_4}) \prod_{i=1}^4 |A_{i,j_i}|^{3/4}$$
(with the convention that $\omega = 0$ when one or more of the four sets is empty).
Since $\omega(A_1,A_2,A_3,A_4) \geq 1/K$, we conclude from the pigeonhole principle that
there exist $j_1,j_2,j_3,j_4 \in \F_2$ with $j_1+j_2+j_3+j_4=0$ such that
$$  \omega(A_{1,j_1},A_{2,j_2},A_{3,j_3},A_{4,j_4}) \prod_{i=1}^4 |A_{i,j_i}|^{3/4}
\geq \frac{1}{K - \frac{1}{10^4}} \big(\prod_{i=1}^4 |A_{i,j_i}|^{3/4}
+ \frac{1}{10^4 K} \prod_{i=1}^4 |A_i|^{3/4}\big).$$
Since $\omega$ is bounded above by $1$, this already implies that $|A_{i,j_i}| \gg K^{-10} |A_i|$ for $i=1,2,3,4$.  We also see that
$$ \omega(A_{1,j_1},A_{2,j_2},A_{3,j_3},A_{4,j_4}) \geq \frac{1}{K - \frac{1}{10^4}}, $$
and the claim follows.
\end{proof}

We may iterate this lemma at most $O(K)$ times to obtain

\begin{corollary}[Large coherently flat subsets] 
Suppose that $A_1,A_2,A_3,A_4 \subseteq \F_2^n$ are non-empty sets with $\omega(A_1,A_2,A_3,A_4) \geq 1/K$ for some $K \geq 1$.  Then there exist $A'_i \subseteq A_i$ with $|A'_i|/|A_i| \geq K^{-O(K)}$, $\omega(A'_1,A'_2,A'_3,A'_4) \geq 1/K$, and such that $(A'_1,A'_2,A'_3,A'_4)$ is coherently $\frac{1}{\sqrt{2K}}$-flat.
\end{corollary}

Theorem \ref{main4} then follows from this corollary and Proposition \ref{well} by repeating the arguments of the previous section.

\end{document}